\documentclass[a4paper,12pt,twoside]{article}
\usepackage{amsthm,amssymb,amsmath,amscd}
\theoremstyle{definition}

\newcommand{\abs}[1]{\left\vert#1\right\vert}

\newcommand{\eps}{\varepsilon}

\newcommand{\Aset}{\mathbb{A}}

\newcommand{\M}{\underline{\text{M}}}
\newcommand{\Mgl}{{(\underline{\text{M}}G)}_{\text{lim}}}
\newcommand{\rc}[1]{C^{*}_{r}({#1})}
\newcommand{\Kt}[1]{K_{j}^{\textrm{top}}({#1})}
\newcommand{\Ku}[1]{K_{i}^{\textrm{top}}({#1})}
\newcommand{\Kv}[1]{K_{i+j}^{\textrm{top}}({#1})}
\newcommand{\Zset}{\mathbb{Z}}

\newcommand{\Cset}{\mathbb{C}}

\newcommand{\Nset}{\mathbb{N}}
\renewcommand{\d}{\textrm{d}}
\begin{document}
\noindent{\footnotesize Mathematical Analysis/\textit{Analyse
Math\'{e}matique}\\ (Topology/\textit{Topologie)}}

\vspace*{1.5cm}

{\LARGE{\bfseries{\noindent A proof of the Baum--Connes
conjecture for\\
reductive adelic groups}}} \\[15pt]
{\large{\bf Paul BAUM$^{a}$, Stephen MILLINGTON$^{b}$, and Roger PLYMEN$^{c}$}} \\[6pt]
$~^{a}$ Department of Mathematics, Pennsylvania State University, University Park, PA16802, USA.\\
$~^{~}$ E-mail:baum@math.psu.edu\\
$~^{b}$ Department of Mathematics, University of Manchester, Manchester, M13 9PL, UK.\\
$~^{~}$ E-mail: stem@quick.freeserve.co.uk\\
$~^{c}$ Department of Mathematics, University of Manchester, Manchester, M13 9PL, UK.\\
$~^{~}$ E-mail: roger@ma.man.ac.uk\\\\
\hspace*{2cm}\hrulefill\hspace*{1.5cm}
\begin{list}{}{\leftmargin 2cm\rightmargin 1.5cm\labelwidth 2cm}
\item[{\bf Abstract.  }] {\small Let $F$ be a global field, $A$ its ring of adeles, $G$
a reductive group over $F$.   We prove the Baum--Connes conjecture
for the adelic group $G(\Aset)$.}

\item{\bf \emph{Une
d\'{e}monstration de la conjecture de Baum--Connes pour les
groupes r\'{e}ductifs ad\`{e}liques}}

\item[{\bf R\'{e}sum\'{e}.\ \ \,}] {\small{\it Soit $F$ un corps global, $A$ son
anneau d'ad\`{e}les, $G$ un groupe r\'{e}ductif sur $F$.   Nous
d\'{e}montrons la conjecture de Baum--Connes pour le groupe ad\`{e}lique $G(\Aset)$.}}
\end{list}
\hspace*{2cm}\hrulefill\hspace*{1.5cm}
\\\\
\noindent{\bf Version Fran\c{c}aise Abr\'{e}g\'{e}e} \\
\\
Soit $X_{1}\subset X_{2}\subset X_{3}\subset \dots$
une suite croissante d'espaces topologiques.   Nous pouvons donner
sur la r\'{e}union $X=\cup X_{n}$ la topologie limite inductive:
un ensemble est ouvert dans $X$ si et seulement si
son intersection avec chaque $X_n$ est ouverte.
Si chaque $X_{n}$ est un espace $T_{1}$ et $E$ est une partie compacte de
$X$, il existe $n$ tel que $E \subset X_{n}$.

Soit $G_1 \subset G_2 \subset G_3 \subset \ldots$
une suite croissante de groupes localement compacts, de base d\'{e}nombrable,
s\'{e}par\'{e}s tel que $G_{n}$ is un sous-groupe ouvert
de $G_{m}$ pour $n\leq m$.  Soit $G=\cup G_{n}$ dans la topologie limite
inductive.    Alors $G$ est un groupe localement compact, de base
d\'{e}nombrable, s\'{e}par\'{e}.   Une base pour la topologie
sur $G$ est donn\'{e}e par l'ensemble des parties ouvertes
dans chaque $G_{n}$.  Si un groupe $G$ est \'{e}gale \`{a}
la r\'{e}union d'une suite croissante de sous-groupes ouverts
alors la topologie sur $G$ est la m\^{e}me que la topologie limite inductive
donn\'{e}e par ces sous-groupes.

Dans \cite{KS} Kasparov et Skandalis construisent un exemple universel
localement compact.   Plus pr\'{e}cisement, soit $\M G$ l'ensemble
des mesures
positives \`{a} support compact sur $G$ avec mesure totale dans l'intervalle
$(1/2,1]$ et muni de la topologie faible *.   A l'aide de l'espace $\M G$
nous construisons un espace pour lequel il est facile
de calculer la $K$-th\'{e}orie topologique d'une r\'{e}union d'une suite
croissante de sous-groupes ouverts.   $\M G$ est un espace m\'{e}trisable
\`{a} base d\'{e}nombrable avec espace quotient m\'{e}trisable
${\M}G/G$.

Si $G$ est la r\'{e}union d'une suite croissante de sous-groupes ouverts
$G_{n}$ alors il existe une suite croissante ferm\'{e}e $\M
G_{1}\subset \M G_{2}\subset \dots$.   L'ensemble $\M G$ est
\'{e}gale \`{a} la r\'{e}union $\cup \M G_{n}$.
Maintenant nous donnons \`{a} $\M G$ une nouvelle topologie, limite inductive
des sous-espaces $\M G_n$.   L'action de $G$ reste continue avec cette nouvelle
topologie.   Remarquons que cet espace n'est pas toujours
localement compact ou m\'{e}trisable --- il tombe en dehors
de la classe des espaces propres au sens de \cite{BCH}.
Il est possible de donner une nouvelle d\'{e}finition de propre ---
plus faible que la d\'{e}finition dans
\cite{BCH}--- pour laquelle cet espace est un exemple universel
pour les actions propres (pour la nouvelle d\'{e}finition) de $G$.
Dans cette Note, il suffit de d\'{e}montrer le r\'{e}sultat suivant.\\

{\sc Th\'{e}or\`{e}me 1.}  L'espace $\M G$, muni de la topologie
inductive limite,
est equivalent par $G$-homotopie \`{a} $\M G$.
\\

Nous d\'{e}montrons aussi les r\'{e}sultats suivants. \\

{\sc Th\'{e}or\`{e}me 2.}   Soit $G$ la r\'{e}union des sous-groupes
$G_n$ localement compacts, de base d\'{e}nombrable,
s\'{e}par\'{e}s, ouverts
et tels que la conjecture de Baum-Connes est vraie pour chaque $G_n$.
Alors la conjecture de Baum-Connes est vraie pour $G$.
\\

Soit $F$ un corps global, c'est-\`{a}-dire ou bien une extension finie
de $\mathbb{Q}$, ou bien une extension finie du corps de fonctions
$\mathbb{F}_p(T)$.
\\

{\sc Th\'{e}or\`{e}me 3.}   Soit $F$ un corps global, $\Aset$ son anneau
d'ad\`{e}les, $G$ un groupe r\'{e}ductif sur $F$.  Alors la
conjecture de Baum-Connes est vraie
pour le groupe ad\`{e}lique $G(\Aset)$.
\\

\noindent{\bf 1. Introduction}

\medskip

Suppose we have an ascending sequence of topological spaces
\[X_{1}\subset X_{2}\subset X_{3}\subset \dots\]
Then we can give the union $X=\cup X_{n}$ the direct limit
topology: a set is open in $X$ if and only if it has open
intersection with each $X_{n}$. If each $X_{n}$ is a $T_{1}$ space
then any compact subset of $X$ lies entirely within some $X_{n}$.
\\
Suppose now that for each $n \in \Nset$ we have a locally compact,
second countable and Hausdorff topological group $G_{n}$, such
that $G_{n}$ is an open subgroup of $G_{m}$ for $n\leq m$. Let
$G=\cup G_{n}$ and furnish this with the direct limit topology.
Then $G$ is a locally compact, second countable, Hausdorff group
in an obvious way. A basis for the topology on $G$ is given by
the collection of open sets in each $G_{n}$. Furthermore if a
group $G$ is equal to the union of an ascending sequence of open
subgroups then the topology on $G$ is the same as the direct limit
topology with respect to these subgroups; note this is not
necessarily the case if the subgroups are not open.

We construct a topological space which is a direct limit and is
homotopic to a universal example for $G$. Using this space we may
express the topological $K$-theory of $G$ as a direct limit. The
adelic groups fit into this framework and the problem reduces to
proving the Baum-Connes conjecture for certain finite products.

Throughout we shall assume that all groups are locally compact,
second countable and Hausdorff.

We would like to thank Paula Cohen, Georges Skandalis and Vincent Lafforgue
for several valuable conversations.    Paul Baum was
partially supported by an NSF grant.

\medskip

\noindent{\bf 2. The space $\M G$}

\medskip

In \cite{KS} Kasparov and Skandalis construct a locally compact
universal example. Namely let $\M G$ be the set of all compactly
supported positive measures on $G$ with total measure in the
interval $(1/2,1]$, topologized with the weak* topology. We shall
use $\M G$ to construct a space for which it is easy to calculate
the topological $K$-theory of an open ascending union of groups.

First we collect some useful properties of $\M G$.
$\M G$ is second countable. This follows from
the fact that $C_{c}(G)$ is separable, which may be seen using the
Stone-Weirstrauss approximation theorem together with the fact
that $G$ is second countable, locally compact and Hausdorff.
So using the Urysohn metrization theorem $\M G$ is metrizable.

Suppose $K\subset \M G$ is compact. We shall show $GK$ is closed.
Let $g_{i}k_{i}$ be a convergent sequence in $GK$ with limit $y$.
We show that $y\in GK$. Let $U$ be a
neighbourhood of $y$ with compact closure.

Since the action is proper (by \cite{CEM} there is no ambiguity
in our use of proper as \M G is locally compact) we have
$\{g\in G:gK\cap \bar{U}\neq\emptyset\}\quad \textrm{is
compact}$.  But because $g_{i}k_{i}$ is convergent this compact
set must contain an infinite number of $g_{i}$ and so $g_{i}$ has
a limit point $g$. Similarly $k_{i}$ all lie in the compact $K$ so
must have a limit point $k$. It is now clear that $y=gk$ as
required.

Using the fact that $\M G$ is normal it is now easy to show that
the quotient space is Hausdorff and locally compact. Also from
the above it is clearly second countable. So using the Urysohn
metrization theorem we can conclude that the quotient space is
metrizable.

\medskip

{\sc Lemma 1.}\quad {\it If $H$ is an open subgroup of
$G$ then $\M H$ is a closed subspace of $\M G$.}

\smallskip

{\sc Proof}\quad A basis for topology is given by the sets
\[(\mu; f_{1},\dots ,f_{n},\eps)=\biggr\{\lambda \in \M G:\abs{\int f_{k}
\:\d\mu - \int f_{k}\:\d\lambda} < \eps,\; \text{for }k=1,\dots ,
n\biggl\}\] where $\mu\in\M G$, $f_{i}\in C_{c}(G)$, and $\eps>0$.
\\
Clearly $\M H\subset \M G$. It is also clear that the open basis
sets $(\lambda; f_{1},\dots,f_{k},\eps)$ of $\M H$ are simply
\[
\biggr\{\lambda \in \M G:\abs{\int f^{G}_{k} \:\d\mu - \int
f^{G}_{k}\:\d\lambda} < \eps,\; \text{for }k=1,\dots ,
n\biggl\}\cap \M H
\]
where each $f^{G}_{k}$ is defined to be $f_{k}$ extended to zero
on $G-H$; giving a continuous function as $H$ is clopen in $G$.
So $\M H$ is a topological subspace of $\M G$.

Now take any $\lambda\not\in \M H$. Then there is some compact
$K\subset G-H$ with $\lambda(K)=\eps>0$. Then the open set
$(\lambda, f_{K}, \eps)$ is a neighbourhood of $\lambda$ which
does not contain any element of $\M H$, where $f_{K}$ is a
compactly supported function which is 1 on $K$ and zero on $H$.
So we see that $\M G-\M H$ is open and hence $\M H$ is closed.
\\

If $G$ is the union of an ascending sequence of open subgroups
$G_{n}$ then there is a closed ascending sequence $\M
G_{1}\subset \M G_{2}\subset \dots$. If we take the union $\cup
\M G_{n}$ it is easy to see that as a set this is the same as $\M
G$. Using the direct limit topology we can then think of $\M G$ as
being retopologized with respect to this ascending sequence. The
action of $G$ remains continuous under this new topology. Note
that this space may not be locally compact or metrizable and so
falls outside the class of spaces we admit as proper in the
sense of \cite{BCH}. However it is possible to give a new
definition of proper---which is weaker than that in
\cite{BCH}---in which this space may be shown to be a universal
example for proper (in this new sense) actions of $G$. However
for our purposes it suffices to prove the following.

\medskip

{\sc Theorem 1.}\quad {\it $\M G$ retopologized in this fashion is
$G$--homotopy equivalent to $\M G$ in the weak* topology.}

\smallskip

{\sc Proof}\quad We write $\Mgl$ to denote $\M G$ with the direct
limit topology. The direct limit topology is finer than the
original topology: this is because each $\M G_{n}$ is closed. So
the identity map $\Mgl\rightarrow \M G$ is obviously  continuous
and $G$-equivariant. We must construct a map in the other
direction.

Take any point $\mu\in \M G$ then as in \cite[1.3]{BCH}
there is a triple
$(U_{\mu}, H_{\mu}, \rho_{\mu})$, where $U_{\mu}$ is a
$G$-invariant open neighbourhood of $\mu$ (in the weak* topology),
$H_{\mu}$ is a compact subgroup of $G$, and $\rho_{\mu}$ is a
$G$-equivariant map from $U_{\mu}$ to $G/H_{\mu}$. As $\M G$ is a
universal example there is a point $\lambda\in\M G$ which is
fixed by $H_{\mu}$. Because the action remains continuous when we
change topology the map $\sigma_{\mu}:G/H_{\mu}\rightarrow \Mgl$
defined by $\sigma_{\mu}(gH_{\mu})=g\lambda$ is well defined and
$G$-equivariant. Define
$\Phi_{\mu}=\sigma_{\mu}\circ\rho_{\mu}:U_{\mu}\rightarrow \Mgl$.
Note that this is a continuous (and $G$-equivariant) map from $U$
with the subspace topology obtained from the weak* topology.

We may cover $\M G$ with such open neighbourhoods and thus get a
cover of ${\M} G/G$ by open sets $\pi(U_{\mu})$. As the quotient
space is paracompact let $\Psi_{\mu}: {\M} G/G \rightarrow [0,1]$ be a
locally finite partition of unity subordinate to this cover.

The map we require is then given by
\[\Xi:\M G\rightarrow\Mgl\quad,\quad\Xi(\nu)=\sum_{\mu}\Psi_{\mu}(\nu)\Phi_{\mu}(\nu)\]
The argument is then completed by noting that any two $G$-maps
from $\M G$ to itself in either topology are $G$-homotopy
equivalent by taking their weighted sums on the interval [0,1].

\medskip

\noindent {\bf 3. K-theory}

\medskip

{\sc Lemma 2.}\quad {\it If $G$ is the union of an ascending
sequence of open subgroups $G_{n}$ then}
\[K_{j}(\rc{G})=\varinjlim K_{j}(\rc{G_{n}})\]

\smallskip

{\sc Proof}\quad We show that $\rc{G}=\varinjlim \rc{G_{n}}$. We
may include $C_{c}(G_{n}) \hookrightarrow C_{c}(G)$ for any $n$
by extending functions to zero. This gives a continuous function
as $G_{n}$ is clopen. Taking $\eps
>0$ and $x\in \rc{G}$ we can get $f\in C_{c}(G)$ which is
$\eps$ close to $x$. But any compact set in $G$ must lie totally
within some $G_{n}$. Hence $f\in C_{c}(G_{n})$
for some $n$ and so $\cup_{n}\rc{G_{n}}$ is dense in $\rc{G}$.
We have

\begin{eqnarray*}
\rc{G}=\varinjlim \rc{G_{n}} \Rightarrow K_{j}(\rc{G})=\varinjlim
K_{j}(\rc{G_{n}})
\end{eqnarray*}
{\sc Lemma 3.}\quad {\it If $G$ is the union of an ascending
sequnce of open subgroups $G_{n}$ then}
\[\Kt{G}=\varinjlim \Kt{G_{n}}\]

\smallskip

{\sc Proof}\quad We have shown that $\Mgl$ is $G$-homotopic to
the universal example $\M G$:
\[\Kt{G}=K_{j}^{G}(\M G)=K_{j}^{G}(\Mgl)\]
By definition 3.13 in \cite{BCH} we have
\[K_{j}^{G}(\Mgl)\ \ =\varinjlim_{\substack{\text{G-compact}\\
Z\subset \M G}}K_{j}^{G}(Z).\] Any $Z\subset \Mgl$ is
$G$-compact if and only if it is the $G$-saturation of a compact
set. But any compact subset of $\Mgl$ lies within some $\M
G_{n}$. So $Z\subset G\cdot \M G_{n}$, for some $n$. We get
\begin{eqnarray*}
K_{j}^{G}(\M G)\ \ =\ \
\varinjlim_{n}\varinjlim_{\substack{\text{G-compact}\\
Z\subset G\cdot \M G_{n}}} K_{j}^{G}(Z) & = &
\varinjlim_{n}K_{j}^{G}(G\cdot \M G_{n}) \\
& = & \varinjlim_{n}K_{j}^{G}(G\times_{G_{n}} \M G_{n})
\end{eqnarray*}

We now appeal to Proposition 5.14 in \cite{CE} which tells us
\[K_{j}^{G}(G\times_{G_{n}} \M G_{n})=K_{j}^{G_{n}}(\M G_{n})=\Kt{G_{n}}\]
as the action of $G_{n}$ on $\M G_{n}$ is proper.

\medskip

{\sc Theorem 2.}  Let $G$ be the union of open subgroups
$G_{n}$ such that the Baum-Connes conjecture is true for each $G_{n}$.
Then the Baum-Connes conjecture is true for $G$.

\smallskip

{\sc Proof}\quad Taking $m\leq n$ we check that the following diagram
commutes
\[\CD
  \Kt{G_{m}} @>\mu_{G_{m}} >> K_{j}(\rc{G_{m}}) \\
  @V VV @V VV \\
  \Kt{G_{n}} @>\mu_{G_{n}} >> K_{j}(\rc{G_{n}})
  \endCD
\]
and use the fact that, for any $G_{m}$-compact $X$, the map $\mu$
factorizes as follows:
\[\CD
  KK_{G_{m}}^{j}(C_{0}(X), \Cset) @>j_{G_{m}} >>
  KK^{j}(C_{0}(X)\rtimes G_{m}, \rc{G_{m}}) @>\#\textbf{1} >>
  KK^{j}(\Cset, \rc{G_{m}})
  \endCD
\]
As the $\mu$ map is defined by direct limits we get the result by
using Lemmas 2 and 3.

\medskip

\noindent{\bf 4. Reductive adelic groups}

\medskip

Let $\mathbb{F}_p(T)$ denote the field of rational functions in the
indeterminate $T$ with coefficients in $\mathbb{F}_p$.   Now
let $F$ be a global field (an $\bf{A}$-field in the sense of
Weil \cite[p.41]{We}).   The field $F$ is a finite algebraic extension
of $\mathbb{Q}$ or a finite algebraic extension of $\mathbb{F}_p(T)$.
We denote by $\Aset$ its ring of adeles, the restricted product of
all the completions $F_{v}$, as in \cite[p.59]{We}.
The field $F$ is a discrete cocompact subfield of the non-discrete
locally compact semisimple commutative ring $\Aset$.
Let $G$ be a reductive group over $F$ and let $G(\Aset)$ denote the group
of adelic points in the algebraic group $G$.
The field $F$ has at most a finite number of infinite places;
it has at least one if it is of characteristic $0$, and none
otherwise.   A place of a global field of characteristic $0$ is
infinite if and only if it lies above the place $\infty$ of $\mathbb{Q}$,
see \cite[p.45]{We}.   Let $S_{\infty}$ be the finite set of
infinite places of $F$ and let $F_{\infty}=\prod_{v\in
S_{\infty}}F_{v}$.
If $v$ is an imaginary place, then the complex reductive group
$G(F_v)$ has also the structure of real reductive group.
We will assume that each local group
$G(F_v)$ with $v \in S_{\infty}$ is a connected Lie group.
Then $G(F_{\infty})$ is a connected real reductive
group. Choose an ordering $v_{1}, v_{2}, \dots$ of the finite
places, and let
\[H_{n}=\prod_{v\leq v_{n}}G(F_{v})\times
\prod_{v>v_{n}}G(\mathfrak{o}_{v})\quad ,\quad
G_{n}=G(F_{\infty})\times H_{n}\]

The groups $G_{n}$ form an inductive system and $G(\Aset)=\bigcup
G_{n}$.   The standard locally compact topology on
$G(\Aset)$ is described in \cite[p.293]{B}.    This topology coincides
with the direct limit topology on $G(\Aset)$ following the open
subgroups $G_n$.

\medskip

{\sc Theorem 3.}   Let $F$ be a global field, $\Aset$
its ring of adeles, $G$ a reductive group over $F$.
Then the Baum--Connes conjecture is true for the adelic group
$G(\Aset)$.

\smallskip

{\sc Proof}\quad Consider the commutative diagram
\[\CD
  \Ku{G(F_{\infty})}\otimes_{\Zset}\Kt{H_{n}} @>\mu_{G(F_{\infty})}\otimes_{\Zset}\mu_{H_{n}}>>
  K_i\rc{G(F_{\infty})}\otimes_{\Zset}K_j\rc{H_{n}} \\
  @V  VV @V \alpha VV  \\
  \Kv{G_{n}} @>\mu_{G_{n}}>> K_{i+j}\rc{G_{n}}
\endCD\]
in which each vertical map is an external product.

We know that $\mu_{G(F_{\infty})}$ is an isomorphism by \cite{LC},\cite{W}.
The ``finite'' product $H_{n}$ admits a 4-tuple satisfying the condition
(HC) of Lafforgue \cite[Definition 1.1]{LS}.   Therefore, by Proposition
1.3 in \cite{LS}, for $t \in \mathbb{R}_+$ sufficiently large,
the Banach space
$S_t(H_n)$ is a good completion of $C_c(H_n)$ and is a subalgebra of
$C^*_r(H_n)$ dense and stable under holomorphic functional calculus.
If $v$ is a finite place then the Euclidean building $\beta G(F_v)$
is a model of the universal example for the local group $G(F_v)$.
Now a Euclidean building is a weakly geodesic and strongly bolic
metric space \cite{K},\cite{KS}.   Let $X_n$ denote the
set of vertices in the product building
$\beta G(F_{v_{1}}) \times \dots\times\beta G(F_{v_{n}})$.
Then $X_n$ admits a metric $d_n$
such that $(X_n, d_n)$ is a weakly geodesic,
strongly bolic and uniformly locally finite metric space on which $H_n$ acts
isometrically, continuously and properly.
By the fundamental result of Lafforgue \cite[Theorem 2.1]{L}, the Baum-Connes
conjecture holds for $H_{n}$ and so
$\mu_{H_{n}}$ is an isomorphism. So $\mu_{G(F_{\infty})}\otimes
\mu_{H_{n}}$ is an isomorphism.

Let $i$ be the mod $2$ dimension of the symmetric space of
$G(F_{\infty})$.   Then $K_{*}\rc{G(F_{\infty})}$ is
concentrated in degree $i$, and is a free abelian group \cite{BC}, \cite{LC},
\cite {W}.  So, by the Kunneth theorem \cite[23.1.3]{BL},
$\alpha$ is an isomorphism. Therefore $\mu_{G_{n}}$ is surjective.

But $G_{n}$ has a $\gamma$ element namely the
external product $\gamma_{G(F_{\infty})}\#\gamma_{H_{n}}$, where
$\gamma_{H_{n}}=\gamma_{v_1}\#\dots\#\gamma_{v_n}$. The individual
$\gamma$-elements $\gamma_{v_1}, \gamma_{v_2}, \dots,
\gamma_{v_n}$ were constructed in \cite{KS}. If $v$ is a finite
place of $F$ then the affine building $\beta G(F_v)$ is a model of the
universal example $\underline{E}G(F_v)$.   The $\gamma$- element
$\gamma_v$ for $G(F_v)$ was constructed in \cite[p.310]{KS0}.   The
$\gamma$-element $\gamma_v$ is an element in the Kasparov ring
$R(G_v) = KK_{G_v}(\Cset, \Cset)$. This implies
the injectivity of $\mu_{G_{n}}$.

Therefore $\mu_{G_{n}}$ is an isomorphism. So the Baum-Connes conjecture
is true for each
open subgroup $G_{n}$. Now apply Theorem 2 to the adelic group
$G(\Aset)=\varinjlim G_{n}$.

\end{document}